\documentclass[preprint,12pt, 3p]{elsarticle}

\usepackage{amssymb}
\usepackage{amsmath}
\usepackage{amsthm}

\usepackage{mathtools}
\usepackage{subcaption}
\usepackage{booktabs}
\newtheorem{thm}{Theorem}
\newtheorem{lem}[thm]{Lemma}
\newtheorem{pro}[thm]{Proposition}

\theoremstyle{definition}

\newtheorem{rmk}[thm]{Remark}
\newproof{pf}{Proof}
\newdefinition{defn}{Definition}

\renewcommand{\eqref}[1]{(\ref{#1})}

\usepackage{xcolor}

\journal{Linear Algebra and its Applications}

\begin{document}

\begin{frontmatter}



\title{Extremal properties and bounds for the generalized algebraic connectivity of graphs in Euclidean spaces\tnoteref{t1, t2}}
\tnotetext[t1]{
© 2025. This manuscript version is made available under the CC-BY-NC-ND 4.0 license.
This is the accepted manuscript of an article published in Linear Algebra and its Applications 727 (2025), 24–36. The final published version is available at https://doi.org/10.1016/j.laa.2025.07.028.
}
\tnotetext[t2]{
This work was partially supported by Universidad de Buenos Aires [grant number UBACyT 20020220100053BA].}

\author[inst1,inst2]{J. Ignacio Alvarez-Hamelin} 
\ead{ihameli@fi.uba.ar}


\affiliation[inst1]{organization={Universidad de Buenos Aires, Facultad de Ingeniería},
            addressline={Av. Paseo Colón 850}, 
            postcode={C1063ACV}, 
            state={Ciudad Autónoma de Buenos Aires},
            country={Argentina}}
            
\author[inst3,inst4]{Juan I. Giribet}
\ead{jgiribet@conicet.gov.ar}
\author[inst3,inst4]{Ignacio Mas}
\ead{imas@udesa.edu.ar}
\author[inst1,inst2]{J. Francisco Presenza\corref{cor1}}
\ead{jpresenza@fi.uba.ar}
\cortext[cor1]{Corresponding author.}

\affiliation[inst2]{organization={CONICET-Universidad de Buenos Aires, INTECIN},
            addressline={Av. Paseo Colón 850}, 
            postcode={C1063ACV}, 
            state={Ciudad Autónoma de Buenos Aires},
            country={Argentina}}

\affiliation[inst3]{organization={Universidad de San Andrés},
            addressline={Vito Dumas 284}, 
            city={Victoria},
            postcode={B1644BID}, 
            state={Provincia Buenos Aires},
            country={Argentina}}

\affiliation[inst4]{organization={CONICET},
            addressline={Godoy Cruz 2290}, 
            postcode={C1425FQB}, 
            state={Ciudad Autónoma de Buenos Aires},
            country={Argentina}}

\begin{abstract}
This article presents new contributions to the study of graph rigidity and its interplay with fundamental graph invariants.
Recently, a quantitative measure of graph rigidity in $\mathbb{R}^d$, termed the \textit{generalized algebraic connectivity}, was introduced.
This development extends the notion of \textit{algebraic connectivity}---the second-smallest eigenvalue of the Laplacian matrix---which is commonly used to quantify graph connectivity.
In this work, we show that the generalized algebraic connectivity is bounded above by the algebraic connectivity. To capture this relationship, we introduce the $d$-\textit{rigidity ratio}, a normalized metric of a graph's rigidity relative to its connectivity.
We also investigate the relationship between rigidity and the \textit{diameter}---a measure of the graph’s overall extent.
In this context, we provide the maximal diameter achievable by rigid graphs and show that \textit{generalized path graphs} serve as extremal examples.
Moreover, we establish a new upper bound for the algebraic connectivity that depends inversely on the diameter and the \textit{vertex connectivity}, improving upon previous bounds.
Finally, we derive an upper bound for the algebraic connectivity of generalized path graphs that asymptotically improves existing ones by a factor of four.
\end{abstract}




\end{frontmatter}


\section{Introduction}
\label{sec:introduction}

Consider an undirected graph $G = (V, E)$ with vertex set $V$ and edge set $E$.
Let $p = (p_i)_{i \in V} \in \mathbb{R}^{d|V|}$ define a realization of $G$ via the mapping $i \in V \textcolor{black}{\mapsto} p_i \in \mathbb{R}^d$.
A pair $(G, p)$ is referred to as a (bar-joint) framework, in which the graph's edges model bars with fixed lengths, and the vertices serve as universal joints.
A fundamental problem in applied geometry is determining whether a framework $(G, p)$ is \textit{rigid}, meaning that rigid motions (compositions of translations and rotations) are the only continuous motions of the vertices that preserve the edge lengths.
Rigidity plays a central role in various applications, as it captures the structural stability of complex networks.
It can be described using the function $f_G: \mathbb{R}^{d|V|} \to \mathbb{R}^{|E|}$ which maps realizations to edge lengths, $f_G(p) = (\Vert p_i - p_j \Vert)_{\{i, j\} \in E}$.
Formally, a framework $(G, p)$ is called rigid if every realization $q$ within some neighborhood of $p$ that satisfies $f_G(q) = f_G(p)$, also has $f_{K_{|V|}}(q) = f_{K_{|V|}}(p)$, where $K_{|V|}$ denotes the complete graph on $V$.
As is common in the literature, we assume that $|V| \geq d + 1$, since otherwise, the framework can be embedded in a proper subspace of $\mathbb{R}^d$, reducing the analysis to a lower-dimensional setting.

The historical applications of this subject are vast, dating back to at least the 19th century with the work of  \citet{Maxwell1864} in structural engineering.
Later, in the 20th century, it found applications in satellite ranging (\citet{Killian1969}), molecular conformation (\citet{Crippen1988}, \citet{Whiteley2005}) and rigidity percolation (\citet{Jacobs1997}).
More recently, the scope of applications has broadened to encompass the localization of sensor networks (\citet{Aspnes2006}) and formation shape control of multi-robot systems (\citet{Anderson2008}, \citet{Krick2009}).

The earliest known rigorous mathematical treatment of rigidity is due to \citet{PollaczekGeiringer1927}, who anticipated the combinatorial characterization of rigidity in $\mathbb{R}^2$ that would later be formalized by \citet{Laman1970}.
The combinatorial characterization of rigidity in $\mathbb{R}^d$, $d \geq 3$, is a major open problem.

On the other hand, \citet{Asimow1978, Asimow1979} established necessary and sufficient algebraic conditions for the rigidity of generic frameworks in arbitrary dimensions.
A framework $(G, p)$ is called \textit{generic} if $p$ maximizes the rank of the \textit{rigidity matrix} among all realizations of $G$.

The rigidity matrix $\mathbf{R}(G, p) \in \mathbb{R}^{|E| \times d|V|}$ has one row for each edge $\{i, j\} \in E$ with the structure $(e_i - e_j)^\mathsf{T} \otimes b_{ij}^\mathsf{T}$, where $\{e_i\}_{i \in V}$ is the standard basis of $\mathbb{
R}^{|V|}$, $\otimes$ denotes the Kronecker product, and  
\begin{equation}
    b_{ij} = 
    \begin{cases}
        (p_i - p_j) / \Vert p_i - p_j \Vert, \quad & p_i \neq p_j, \\
        \boldsymbol{0}_d, \quad & \text{otherwise}.
    \end{cases}
    \label{eq:bearing_vector}
\end{equation}
The null space of the rigidity matrix, $\mathcal{N}(\mathbf{R}(G, p))$, contains the subspace of \textit{infinitesimal rigid motions}, defined as
\begin{equation*}
    \mathcal{T}(p) = \{(u_i)_{i \in V} \in \mathbb{R}^{d|V|} : u_i = A p_i + q, \ (A, q) \in \mathfrak{so}(d) \times \mathbb{R}^d, \, \text{for all} \ i \in V\}.
\end{equation*}
A framework $(G, p)$ is called \textit{infinitesimally rigid} if $\mathcal{T}(p) = \mathcal{N}(\mathbf{R}(G, p))$\textcolor{black}{; or, equivalently}, if $\mathrm{rank}\left(\mathbf{R}(G, p)\right) = d|V| - \binom{d+1}{2}$, which follows since the dimension of $\mathcal{T}(p)$ is $\binom{d+1}{2}$, provided that $\{p_i\}_{i \in V}$ affinely spans $\mathbb{R}^d$.
It is well known that infinitesimal rigidity is a stronger property than rigidity, and that both notions coincide for generic frameworks.

Until recently, framework rigidity was posed as a yes or no question.
Motivated by the deployment of sensor networks and the control of multi-robot systems, quantitative measures of framework rigidity were developed.
This arises, for instance, from the need to assess robustness and reliability in estimating or controlling vertex positions based on range measurements.
To this end, \citet{Zhu2009} and \citet{Zelazo2015} introduced the positive semi-definite \textit{stiffness matrix} $\mathbf{S}(G, p) \coloneqq \mathbf{R}(G, p)^\mathsf{T} \mathbf{R}(G, p)$ and showed that $(G, p)$ is infinitesimally rigid if and only if $\lambda_{\binom{d+1}{2} + 1}(\mathbf{S}(G, p)) > 0$, where $\lambda_k$ denotes the $k$-th smallest eigenvalue of $\mathbf{S}(G, p)$.
The latter is known as the \textit{rigidity eigenvalue} of $(G, p)$, which acts as an indicator of its level of rigidity.
More recently, \citet{Jordan2022} introduced a measure associated with the rigidity of the underlying graph, termed \textit{generalized algebraic connectivity}.
This is defined for $d \geq 1$ and $|V| \geq d + 1$ as
\begin{equation}
    a_d(G) \coloneqq \sup \Big\{\lambda_{\binom{d+1}{2}+1}(\mathbf{S}(G, p)) : p \in \mathbb{R}^{d|V|}\Big\}.
    \label{eq:d_algebraic_connectivity}
\end{equation}
It follows that $a_d(G) > 0$ if and only if $G$ admits a realization $p \in \mathbb{R}^{d|V|}$ such that $(G, p)$ is infinitesimally rigid.
In this case, $(G, q)$ is infinitesimally rigid for every $q$ in a dense open subset of $\mathbb{R}^{d|V|}$ whose complement has Lebesgue measure zero (see \citet{Asimow1978}), and $G$ is said to be \textit{generically rigid} in $\mathbb{R}^d$.
For $d=1$, $a_1(G)$ coincides with the usual algebraic connectivity of graphs $a(G)$, meaning that generic rigidity in $\mathbb{R}^1$ and connectivity are equivalent.

This article presents new results on the generalized algebraic connectivity of graphs in Euclidean spaces, and its relationship with other graph invariants such as the algebraic connectivity, the diameter, and the vertex connectivity.
The structure and main contributions are listed as follows. 
Section \ref{sec:preliminaries} completes the introduction by presenting the notation used and additional definitions.
Section \ref{sec:ratio} presents the $d$-\textit{rigidity ratio}, $\rho_d(G) \coloneqq a_d(G) / a_1(G)$, as a way to normalize a graph's generalized algebraic connectivity with respect to its connectivity. We show that $a_d(G) \leq a_1(G)$ (and hence $\rho_d(G) \leq 1$), \textcolor{black}{and discuss potential extremal examples}.
Section \ref{sec:rigidity_diameter} is devoted to relating rigidity and the graph's diameter.
Firstly, we determine that the diameter of a generically rigid graph in $\mathbb{R}^d$ with $n$ vertices satisfies $\mathrm{diam}(G) \leq \lceil \frac{n - 1}{d}\rceil$.
Moreover, we show that \textit{generalized path graphs} \textcolor{black}{(see \eqref{eq:path_graph})} attain this maximum.
Secondly, we present a new upper bound for the algebraic connectivity of a graph in terms of the diameter and the vertex connectivity, $\kappa(G)$, specifically

\begin{equation*}
    a_1(G) \leq \frac{12 |E(G)|}{\kappa(G) \mathrm{diam}(G) (\mathrm{diam}(G)-1) (\mathrm{diam}(G)-2) + 6 \mathrm{diam}(G)^2}.
\end{equation*}
Finally, the algebraic connectivity of generalized path graphs is studied, and the upper bound $a_1(P_{n, d}) \leq a_1(C_{2n, d})$ is provided\textcolor{black}{, where $C_{n, d}$ is the \textit{generalized cycle graph} (see \eqref{eq:cycle_graph}). This bound }
improves upon \citet[Prop. 7.8]{Lew2025} by a factor of four (asymptotically, as the graph size increases).
Numerical examples comparing the new bounds with similar ones from the literature are provided for the generalized path graph.
Section \ref{sec:conclusion} provides a summary and concluding remarks.

\section{Preliminaries}\label{sec:preliminaries}
\subsection{Notation}

Let $\mathbb{N}$ be the set of positive integers, and $[n] \coloneqq \{1, \ldots, n\}$.
For two vectors $u, v \in \mathbb{R}^d$, $\langle u, v \rangle$ denotes the standard inner product, and $\Vert u \Vert = \langle u, u \rangle^{\frac{1}{2}}$.
The $(d-1)$-sphere is $\mathcal{S}^{d-1} = \{u \in \mathbb{R}^d : \Vert u \Vert = 1\}$.
The $d \times 1$ vectors of all zeros and all ones are symbolized as $\boldsymbol{0}_d$ and $\boldsymbol{1}_d$, respectively.
The $d \times d$ identity matrix is indicated as $\mathbf{I}_d$.
The eigenvalues of an $\mathbb{R}^{n \times n}$ symmetric matrix are arranged as $\lambda_1 \leq \ldots \leq \lambda_n$, counted with multiplicities.
The cardinality of a set $S$ is denoted as $|S|$.
Let $\{u_i\}_{i \in S} \subset \mathbb{R}^d$ where $S$ is a finite set of ordered indices.
We define $(u_i)_{i \in S} = (u_{i_1}, \ldots, u_{i_{|S|}}) \in \mathbb{R}^{d|S|}$ as the column-wise stacking of the vectors $u_i$, preserving their order in $S$.
The space of real skew-symmetric matrices is denoted as $\mathfrak{so}(d)$.

\subsection{Graphs}
We consider finite undirected graphs $G = (V, E)$ with no loops or multiple edges.
For notational simplicity, vertices $V = V(G) = [n]$ are numbered, and edges $E = E(G) \subseteq \{\{i, j\} \subseteq V : i \neq j\}$ are ordered lexicographically.
The degree of a vertex $i \in V$, that is, the number of edges that are incident to it, is denoted as $\mathrm{deg}(i)$.
A path of length $m$ in $G$ is a sequence of distinct vertices $\{i_0, \ldots, i_m\}$ such that $\{i_k, i_{k+1}\} \in E$ for each $k=0, \ldots, m-1$.
A graph is connected if, for every pair $i,j \in V$, there is a path that has them as endpoints. 
The following definitions assume that the graph $G$ is connected.
The distance between two vertices $i, j \in V$, denoted $\delta(i, j)$, is the length of the shortest path connecting them.
The graph's diameter is defined as $\mathrm{diam}(G) = \max_{i, j \in V} \delta(i, j)$.
The following definition classifies connected graphs according to the minimum number of vertices that must be removed to disconnect them.
The subsequent theorem establishes a connection between this property and the number of vertex-disjoint paths between vertices; see ({\citet[Th. 3.4.1]{Godsil2013}}).

\begin{defn}
    Let $G$ be a connected graph with $n$ vertices, and let $k \in \mathbb{N}$ such that $k < n$. Then, $G$ is said to be $k$-connected if it remains connected after the removal of any subset of fewer than $k$ vertices.
    \textit{Vertex connectivity} $\kappa(G)$ is the largest $k$ for which $G$ is $k$-connected. 
    \label{def:vertex_connectivity}
\end{defn}

\begin{thm}[Menger]
    Let $G = (V, E)$ be a connected graph with $n \geq 2$ vertices, and let $i, j \in V$ be distinct.
    Then, the number of paths having $i$ and $j$ as endpoints that are internally vertex-disjoint is no less than the vertex connectivity $\kappa(G)$.
    \label{th:menger}
\end{thm}

Graph connectivity can also be assessed via the $n \times n$ Laplacian matrix (\citet{Godsil2013}), defined such that the $(i, j)$-th entry 
\begin{equation*}
    [\mathbf{L}(G)]_{i, j} =
    \begin{cases}
        -1, \quad & \{i, j\} \in E, \\
        \mathrm{deg}(i) \quad & i = j, \\
        0, \quad & \text{otherwise}.
    \end{cases}
\end{equation*}
The Laplacian is symmetric positive semi-definite, $\lambda_1(\mathbf{L}(G)) = 0$ and $\boldsymbol{1}_n \in \mathcal{N}(\mathbf{L}(G))$.
Moreover, if $n \geq 2$, $G$ is connected if and only if $a(G) \coloneqq \lambda_2(\mathbf{L}(G)) > 0$. 
This eigenvalue is known as the \textit{algebraic connectivity} of $G$, and serves as an indicator of its level of connectivity.
As previously mentioned, $a_1(G) = a(G)$, establishing the equivalence between connectivity and generic rigidity in $\mathbb{R}^{1}$.
The following theorem extends this result by relating such properties in higher dimensions.
It follows directly from \citet[Th. 61.1.4]{Goodman2017}, to which the reader is referred for a proof.

\begin{thm}
    Let $d \in \mathbb{N}$ and let $G$ be a graph with $n \geq d+1$ vertices. If $G$ is generically rigid in $\mathbb{R}^d$, then $\kappa(G) \geq d$.
    \label{th:d_connectivity}
\end{thm}

\section{$d$-Rigidity Ratio}
\label{sec:ratio}

The advantages of having the generalized algebraic connectivity \eqref{eq:d_algebraic_connectivity} extend beyond merely assessing generic rigidity in any dimension.
Considering that rigidity can be viewed as an extension of connectivity to higher dimensions,
we propose the $d$-\textit{rigidity ratio} 
\begin{equation*}
    \rho_d(G) \coloneqq \frac{a_d(G)}{a_1(G)}, \quad d \geq 2
\end{equation*}
as a normalized quantifier of the rigidity level of a connected graph in $\mathbb{R}^d$.

\subsection{The $d$-rigidity ratio lies within the unit interval}

When studying the bounds of the $d$-rigidity ratio for $d \geq 2$, it immediately follows that $\min\{\rho_d(G) : G \text{ is connected}\} = 0$, and the lower bound is attained by all connected graphs that are not generically rigid in $\mathbb{R}^d$.
If $n \geq d+1$, the \textit{path graph} $P_n$ defined by $V(P_n) =[n]$ and $E(P_n) =\{\{i, j\} : |i - j| = 1\}$ serves as such an example.
Two fundamental questions emerge: first, whether $\{\rho_d(G) : G \text{ is connected}\}$ is bounded above; and second, determining the existence of maximizing graphs.
The resolution of the first question is provided by the main result of this section, Theorem \ref{th:algebraic_connectivity_bound}, which shows that $\rho_d(G) \leq 1$. 
The special case $d=2$ was previously proved by \citet[Th. 4.2]{Jordan2022}.
The general case for $d \geq 2$ was first introduced in a preliminary (unpublished) version of this work available at \cite{Presenza2022b}, and was later independently proved by \citet{Lew2025} using probabilistic methods.

\begin{lem}
    Let $d \in \mathbb{N}$ and $(G, p)$ be a $d$-dimensional framework with $n \geq d+1$ vertices. 
    Let $v^\star = (v^\star_i)_{i \in [n]} \in \mathbb{R}^n$ be a unit norm eigenvector associated with $\lambda_2(\mathbf{L}(G))$, and let $w \in \mathcal{S}^{d-1}$. 
    Define $u^\star = v^\star \otimes w$, then
    \begin{equation*}
        \langle \mathbf{S}(G, p) u^\star, u^\star \rangle \leq a_1(G).
    \end{equation*}
    Moreover, equality is achieved if and only if $p$ is such that $b_{ij} = \pm w$ for all $\{i, j\} \in E$, with $b_{ij}$ defined as in \eqref{eq:bearing_vector}.
    \label{lem:quad_form_bound}
\end{lem}
\begin{pf}
    By expanding the quadratic form of the stiffness matrix as a sum (see \citet[Sec. 3]{Jordan2022}), one gets
    \begin{equation*}
        \langle \mathbf{S}(G, p) u^\star, u^\star \rangle = \sum_{\{i, j\} \in E} \langle (v^\star_i - v^\star_j) w, b_{ij} \rangle^2
        \leq \sum_{\{i, j\} \in E} (v^\star_i - v^\star_j)^2 = \langle \mathbf{L}(G) v^\star, v^\star\rangle = \lambda_2(\mathbf{L}(G))
    \end{equation*}
    and equality holds if and only if $|\langle w, b_{ij} \rangle| = 1$ for all $\{i, j\} \in E$.
    \qed
\end{pf}

\begin{thm}
    Let $d  \geq 2 \in \mathbb{N}$ and let $G$ be a connected graph with $n \geq d + 1$ vertices.
    Then,
    \begin{equation*}
        a_d(G) \leq a_1(G).
    \end{equation*}
    \label{th:algebraic_connectivity_bound}
\end{thm}
\begin{pf}
    Consider a generically rigid graph $G$ in $\mathbb{R}^d$, otherwise $a_d(G) = 0$. 
    Let $p \in \mathbb{R}^{dn}$ be a realization for which $(G, p)$ is infinitesimally rigid.
    Let $v^\star = (v^\star_i)_{i \in [n]} \in \mathbb{R}^n$ be a unit norm eigenvector associated with $\lambda_2(\mathbf{L}(G))$.
    Let $\tilde{w} = \sum_{i=1}^n v^\star_i p_i \in \mathbb{R}^d$ and define
    \begin{equation*}
        w = 
        \begin{cases}
            \tilde{w} / \Vert \tilde{w} \Vert, \quad & \Vert \tilde{w} \Vert \neq 0, \\
            \text{any vector in} \ \mathcal{S}^{d-1}, \quad & \Vert \tilde{w} \Vert = 0.
        \end{cases}
    \end{equation*}
    \textcolor{black}{Let $u^\star \coloneqq v^\star \otimes w$}. We show that $u^\star$ is orthogonal to $\mathcal{T}(p)$.
    To see this, consider $x \in \mathcal{T}(p)$; then $x = (A p_i + q)_{i \in [n]}$ where $A \in \mathfrak{so}(d)$ and $q \in \mathbb{R}^d$.
    Thus, 
    \begin{equation*}
        \langle u^\star, x \rangle = \sum_{i=1}^n \langle v_i^\star w, A p_i + q \rangle = 
        \left \langle w, A \left(\sum_{i=1}^n v_i^\star  p_i\right) \right \rangle + \left(\sum_{i=1}^n v_i^\star\right) \langle w, q \rangle =  
        \langle w, A \tilde{w} \rangle =
        0,
    \end{equation*}
    where we used that $\sum_{i=1}^n v_i^\star = 0$ due to the orthogonality between $v^\star$ and $\boldsymbol{1}_n$ (as eigenvectors of $\mathbf{L}(G)$), and that $\langle w, A \tilde{w} \rangle = \langle w, A w \rangle = 0$, since $A$ is skew-symmetric. 
    The claim follows since
    \begin{equation}
        \lambda_{\binom{d+1}{2}+1}(\mathbf{S}(G, p)) = \min\!\big\{\langle \mathbf{S}(G, p) x, x \rangle : \Vert x \Vert\!=\!1, x\!\perp\!\mathcal{T}(p)\big\} \leq \langle \mathbf{S}(G, p) u^\star, u^\star \rangle \leq a_1(G),
        \label{eq:th_ineq_2}
    \end{equation}
    where we used Lemma \ref{lem:quad_form_bound}.
    \qed
\end{pf}

The question of whether there exist graphs that attain the equality $a_d(G) = a_1(G)$ for $d \geq 2$ remains open.
Remarks \ref{rmk:complete_not_sharp} and \ref{rmk:star_not_sharp} indicate that complete graphs \textcolor{black}{$K_n$} and generalized star graphs \textcolor{black}{$S_{n, d}$} (\textcolor{black}{consisting of a clique
of size $d$, and $n-d$ additional vertices, each adjacent to, and only to, all the vertices of the clique, see} \citet[Sec. 1]{Lew2025}) \textcolor{black}{do not achieve this equality}.
This conclusion follows from the existing results by \citet[Th. 4.4]{Jordan2022}, \citet[Ths. 1.2 and 1.6]{Lew2023}, and \citet[Prop. 7.6]{Lew2025}.
\textcolor{black}{
Observe that $K_2$ is not considered, since $a_d(K_2)$ is defined only for $d=1$.
}

\begin{rmk}
    It holds that 
    \begin{align*}
        \rho_2(K_n) = \tfrac{1}{2}, \ n \geq d+1, \
        \rho_d(K_{d+1}) = \tfrac{1}{d + 1}, \ d \geq 3, \ \text{and} \
        \rho_d(K_n) < \tfrac{1}{d-1}, \ d \geq 3, n \geq d+2
    \end{align*}
    \label{rmk:complete_not_sharp}
\end{rmk}
\begin{rmk}
    Let $S_{n, d}$ be the It holds that $\rho_d(S_{n, d}) = \frac{1}{d}$ for all $d \geq 2, n \geq d+2$.
    \label{rmk:star_not_sharp}
\end{rmk}

\section{Impact of Graph Diameter on Rigidity Properties}
\label{sec:rigidity_diameter}
The interplay between graph connectivity and the diameter has been a subject of considerable interest over the years.
In modeling communication networks, the diameter is indicative of message delay throughout the network; thus, it is crucial for both performance analysis and optimization.
For instance, it plays an important role in the performance of distributed consensus-based protocols, see \citet{Hendrickx2014}.
Several results relating graph connectivity and the diameter can be found in the book by \citet{Chung1996}, as well as in the works of \citet{Alon1985}, \citet{Mohar1991}, \citet{Deabreu2007}; however, there are still many open problems.
On the other hand, the relationship between rigidity and the diameter is, to our knowledge, unexplored.
This question can be particularly valuable for the operation of sensor networks and multi-robot systems (\citet{Zelazo2015}, \citet{Amani2020}, and \citet{Presenza2022a}).

\subsection{Maximum Diameter of Rigid Graphs}
\label{sec:max_diam_rigid}
Here we explore the maximum diameter achievable for generically rigid graphs in each dimension.

\begin{pro}
    Let $d \in \mathbb{N}$, and let $G$ be a generically rigid graph in $\mathbb{R}^d$ with $n \geq d + 1$ vertices. Then, 
    \begin{equation}
        \mathrm{diam}(G) \leq \Big\lceil \frac{n-1}{d} \Big\rceil.
        \label{eq:rigid_maximal_diameter}
    \end{equation}
    \label{pro:rigid_maximal_diameter}
\end{pro}
\begin{pf}
    We know from Theorem \ref{th:d_connectivity} that $G$ is $d$-connected.
    Hence, according to \citet[Th. 1]{Caccetta1992}, $G$ has a diameter less than or equal to $\lceil \frac{n-1}{d} \rceil$.
    \qed
\end{pf}

The sharpness of \eqref{eq:rigid_maximal_diameter} is demonstrated by Proposition \ref{pro:path_diameter}, where the extremal example given is the \textit{generalized path graph} $P_{n, d}$, defined for $n \geq d+1$ such that
\begin{equation}
    V(P_{n, d}) = [n] \quad \text{and} \quad 
    E(P_{n, d}) = \{\{i, j\} : |i - j| \in \{1, \ldots, d\}\}.
    \label{eq:path_graph}
\end{equation}
This result generalizes the fact that $P_{n, 1} = P_n$ attains the maximum diameter among all connected graphs with $n$ vertices.

\begin{pro}
    Let $d,n \in \mathbb{N}$ such that $n \geq d+1$, and let $P_{n, d}$ be the generalized path graph \eqref{eq:path_graph}. Then, 
    \begin{equation*}
        \mathrm{diam}(P_{n, d}) = \Big\lceil \frac{n-1}{d} \Big\rceil.
    \end{equation*}
    \label{pro:path_diameter}
\end{pro}
\begin{pf}
    We proceed by induction for $n \geq d+1$.
    For $n = d+1$, we have $\mathrm{diam}(P_{d+1, d}) = 1 = \lceil (n-1) / d \rceil$.
    Now, consider $n \geq d+1$ and assume that $\mathrm{diam}(P_{n, d}) = \lceil (n-1) / d \rceil$.
    Without loss of generality, we label the vertices of $P_{n, d}$ as $V = [n]$. 
    It follows that $\delta(1, i) = \lceil (i-1) / d \rceil$.
    When a new vertex $j = n + 1$ is added to $P_{n, d}$ to obtain $P_{n+1, d}$, it will be connected with the vertices $n, n-1, \ldots, n - d + 1$; therefore 
    \begin{equation*}
        \mathrm{diam}(P_{n+1, d}) = \delta(1, n+1) = \delta(1, n - d + 1) + 1 = \lceil (n-d) / d \rceil + 1 = \lceil n / d \rceil.       
    \end{equation*}
    \qed
\end{pf}

\subsection{Diameter and the Generalized Algebraic Connectivity}
\label{sec:algebraic_connectivity_diam_bound}

Here, we present an upper bound for $a_1(G)$ (and thus for $a_d(G)$) that is inversely related to the diameter and to the vertex connectivity.

\begin{thm}
    Let $G$ be a connected graph with $n \geq 2$ vertices, vertex connectivity $\kappa(G)$ and diameter $\mathrm{diam}(G)$.
    Then,
    \begin{equation}
        a_1(G) \leq \frac{12 |E(G)|}{\kappa(G)  \mathrm{diam}(G) (\mathrm{diam}(G)-1) (\mathrm{diam}(G)-2) + 6 \mathrm{diam}(G)^2}.
        \label{eq:diameter_inequalities_2}
    \end{equation}
    \label{th:diameter_bound}
\end{thm}
\begin{pf}
    In the following, we omit the dependence on $G$ in $E(G)$, $\kappa(G)$ and $\mathrm{diam}(G)$ for simplicity.
    Let $a, b$ be two vertices such that $\delta(a, b) = \mathrm{diam}$.
    Define $v = (v_i)_{i \in [n]}$ such that $v_i = \delta(a, i)/\mathrm{diam}$.
    It follows that $v_a = 0, v_b = 1$ and $0 < v_i \leq 1$ for all $i \notin \{a, b\}$. 
    Also, define $\tilde{v} = (\tilde{v}_i)_{i \in [n]}$ such that $\tilde{v}_i = v_i - \bar{v}$ where $\bar{v} = 1/n \sum_{i=1}^n v_i$.
    It follows that $\langle \boldsymbol{1}_n, \tilde{v} \rangle = 0$ and that
    \begin{equation*}
        \langle \mathbf{L}(G) \tilde{v}, \tilde{v} \rangle = \sum_{\{i, j\} \in E} (\tilde{v}_i - \tilde{v}_j)^2 = \sum_{\{i, j\} \in E} (v_i - v_j)^2 \leq \frac{|E|}{\mathrm{diam}}^2,
    \end{equation*}
    since $\{i, j\} \in E$ implies $(v_i - v_j)^2 \in \{0, \mathrm{diam}^{-2}\}$.
    Now, due to Theorem \ref{th:menger} (Menger), consider $\kappa$ internally vertex-disjoint paths that connect $a$ and $b$, $\Pi_1, \ldots, \Pi_{\kappa}$.
    Each path has length greater than or equal to $\mathrm{diam}$ and contains at least one vertex $i$ such that $\delta(a, i) = \ell$, for each $\ell = 0, \ldots, \mathrm{diam}$.
    Define $\Pi_j^* = \Pi_j \setminus \{a, b\}$, then $\Pi^*_j \cap \Pi^*_{j'} = \emptyset$ if $j \neq j'$. Hence,
    \begin{align*}
        \langle \tilde{v}, \tilde{v} \rangle = \sum_{i=1}^n  (v_i - \bar{v})^2 &\geq (v_a - \bar{v})^2 + (v_b - \bar{v})^2 + \sum_{j=1}^{\kappa} \sum_{i \in \Pi_j^*} (v_i - \bar{v})^2 \\
        &\geq \bar{v}^2 + (1 - \bar{v})^2 + \sum_{j=1}^{\kappa} \sum_{\ell=1}^{\mathrm{diam}} - 1 \left(\frac{\ell}{\mathrm{diam}} - \bar{v}\right)^2 \\
        &\geq \frac{1}{2} + \kappa\sum_{\ell=1}^{\mathrm{diam}-1} \left(\frac{\ell}{\mathrm{diam}} - \bar{v}\right)^2 \\
        & \geq \frac{1}{2} + \kappa\sum_{\ell=1}^{\mathrm{diam}-1} \left(\frac{\ell}{\mathrm{diam}} - \frac{1}{2}\right)^2 \\
        &= 
        \frac{1}{2} + \kappa\frac{(\mathrm{diam}-1)(\mathrm{diam}-2)}{12 \mathrm{diam}}.
    \end{align*}
    The first inequality holds since vertices not contained in any path $\Pi_j$ are removed from the summation.
    The second one, since for each $1 \leq \ell \leq \mathrm{diam}-1$ only one vertex $i$ such that $\delta(a, i) = \ell$ is counted per path.
    The third one is due to the fact that $\bar{v}^2 + (1 - \bar{v})^2 \geq 1/2$.
    And the last one holds since $1/2$ is the mean value of $\{\ell/\mathrm{diam} : \; 1 \leq \ell \leq \mathrm{diam}-1\}$, and the squared deviations from $\bar{v}$ are no less than the squared deviations from the mean.
    
    As the final step,
    \begin{equation*}
        a_1(G) = \min \Bigg\{\frac{ \langle \mathbf{L}(G) x, x \rangle}{\langle x, x \rangle} : \Vert x \Vert = 1, x \perp \boldsymbol{1}_n\Bigg\} \leq \frac{\langle \mathbf{L}(G) \tilde{v}, \tilde{v} \rangle}{\langle \tilde{v}, \tilde{v} \rangle},    
    \end{equation*}
    and the upper bound \eqref{eq:diameter_inequalities_2} holds.
    \qed
\end{pf}

It can be observed that $K_2$ achieves the equality in \eqref{eq:diameter_inequalities_2}, demonstrating its sharpness.
However, it is of interest to find more complex graphs approaching the equality.

\subsubsection{Generalized Path Graphs}
\label{sec:path_graphs}

Following the prior discussion, and the fact that $a_1(G) \geq a_1(P_{n, 1})$ for every connected graph $G$ with $n$ vertices (see \citet{Fiedler1973}), it is natural to conjecture that $a_d(G) \geq a_d(P_{n, d})$ for every rigid graph $G$ in $\mathbb{R}^d$ with $n$ vertices.
Thus, it is of interest to derive precise expressions for $a_d(P_{n, d})$.
In this direction, \citet[Prop. 7.8]{Lew2025} showed that  $a_1(P_{n, d}) \leq a_1(C_{n, d})$ where $C_{n, d}$ is called the \textit{generalized cycle graph}, defined for $n \geq d+1$ such that
    \begin{equation}
    \begin{split}
        V(C_{n, d}) = [n] \quad \text{and} \quad 
        E(C_{n, d}) = \{\{i, j\}: \; |i - j| \in \{1, \ldots, d, n-d, \ldots, n-1\}\}.
    \end{split}
    \label{eq:cycle_graph}
    \end{equation}

Motivated by the identity $a_1(P_{n, 1}) = a_1(C_{2n, 1})$ (see \citet[Lem. 5.6.1]{Spielman2025}), we provide Proposition \ref{pro:path_cycle} which presents a new upper bound for $a_1(P_{n, d})$  considering $C_{2n, d}$ instead of $C_{n, d}$.
Proposition \ref{pro:bounds_ratio} shows that it improves upon \citet[Prop. 7.8]{Lew2025} by a factor of four, asymptotically as $n$ grows.

To do this, define the vectors $u_{(n)}, v_{(n)} \in \mathbb{R}^n$ such that
\begin{equation}
    u_{(n)} = \left(\sqrt{\tfrac{2}{n}} \cos\left(\tfrac{2\pi}{n}\left(i - \tfrac{1}{2}\right)\right)\right)_{i \in [n]} \quad \text{and} \quad 
    v_{(n)} = \left(\sqrt{\tfrac{2}{n}} \cos\left(\tfrac{\pi}{n}\left(i - \tfrac{1}{2}\right)\right)\right)_{i \in [n]}.
    \label{eq:eigenvector_cycle}
\end{equation}

\begin{rmk}
    $\Vert u_{(n)} \Vert = \Vert v_{(n)} \Vert = 1$, $\langle \boldsymbol{1}_n, u_{(n)} \rangle = \langle \boldsymbol{1}_n, v_{(n)} \rangle = 0$ and $u_{(2n)} = \left(\frac{v_{(n)}}{\sqrt{2}}, -\frac{v_{(n)}}{\sqrt{2}}\right)$.
    \label{rmk:eigenvector_cycle_2}
\end{rmk}

\begin{lem}[{\citet[Chap. 3]{Gray2006}}]
    Let $d, n \in \mathbb{N}$ such that $n \geq d+1$ and $C_{n, d}$ be defined as in \eqref{eq:cycle_graph}. Then, $\mathbf{L}(C_{n, d})$ is a circulant matrix; hence,
    \begin{equation*}
        a_1(C_{n, d}) = 
        \begin{cases}
            n, \quad & n \leq 2d+1, \\
            \sum_{k=1}^d 2\left(1 - \cos(\tfrac{2 k \pi}{n})\right), \quad & n \geq 2d+2,
        \end{cases}
    \end{equation*}
    and $u_{(n)}$ \eqref{eq:eigenvector_cycle} is an associated unit eigenvector.
    \label{lem:cycle_spectrum}
\end{lem}

\begin{lem}
    Let $d, n \in \mathbb{N}$ such that $n \geq d+1$.
    Consider the graphs $P_{n, d}$ \eqref{eq:path_graph} and $P^*_{n, d}$, where $V(P^*_{n, d}) = \{n+1, \ldots, 2n\}$ and $\{i, j\} \in E(P^*_{n, d})$ if and only if $\{i-n, j-n\} \in E(P_{n, d})$.
    \textcolor{black}{Then},
    \begin{equation*}
        \langle \mathbf{L}(P_{n, d}) v_{(n)}, v_{(n)} \rangle \leq \langle \mathbf{L}(C_{2n, d}) u_{(2n)}, u_{(2n)} \rangle.
    \end{equation*}
    with $v_{(n)}, u_{(2n)}$ taken from \eqref{eq:eigenvector_cycle}.
    \label{lem:path_isomorphism}
\end{lem}
\begin{pf}
    Observe that\textcolor{black}{, with the usual graph union,} $\mathbf{L}(P_{n, d} \cup P^*_{n, d}) = \mathbf{I}_2 \otimes \mathbf{L}(P_{n, d})$.
    This holds since $\textcolor{black}{V(P_{n, d}) \cap V(P^*_{n, d})} = \emptyset$ and $\mathbf{L}(P^*_{n, d}) = \mathbf{L}(P_{n, d})$. 
    From this observation and from Remark \ref{rmk:eigenvector_cycle_2}, it follows that
    \begin{equation*}
        \langle \mathbf{L}(P_{n, d}) v_{(n)}, v_{(n)} \rangle = \langle \mathbf{L}(P_{n, d} \cup P^*_{n, d}) u_{(2n)}, u_{(2n)} \rangle \leq \langle \mathbf{L}(C_{2n, d}) u_{(2n)}, u_{(2n)} \rangle.
    \end{equation*}
    For the inequality, we used \textcolor{black}{$V(P_{n, d} \cup P^*_{n, d}) = V(C_{2n, d})$} and $E(P_{n, d} \cup P^*_{n, d}) \subset E(C_{2n, d})$ and the monotonicity of the Laplacian quadratic form with respect to the addition of edges.
    \qed
\end{pf}

\begin{pro}
    Let $d, n \in \mathbb{N}$ such that $d \geq 2$ and $n \geq d+2$, and let $P_{n, d}$ and $C_{n, d}$ be defined as in \eqref{eq:path_graph} and \eqref{eq:cycle_graph}, respectively.
    Then,
    \begin{equation}
        a_1(P_{n, d}) \leq a_1(C_{2n, d}).
        \label{eq:path_cycle}
    \end{equation}
    \label{pro:path_cycle}
\end{pro}
\begin{pf}
    Consider $v_{(n)}$ and $u_{(n)}$ from \eqref{eq:eigenvector_cycle}.
    Then,
    \begin{align*}
        a_1(P_{n, d}) &= \min \big\{\langle \mathbf{L}(P_{n, d}) x, x \rangle : \Vert x \Vert = 1, x \perp \boldsymbol{1}_n\big\} \\
        &\leq \langle \mathbf{L}(P_{n, d}) v_{(n)}, v_{(n)} \rangle \leq \langle \mathbf{L}(C_{2n, d}) u_{(2n)}, u_{(2n)} \rangle = a_1(C_{2n, d}),
    \end{align*}
    where we used Lemmas \ref{lem:cycle_spectrum} and \ref{lem:path_isomorphism}.
    \qed
\end{pf}

\begin{pro}
    Let $d \geq 2 \in \mathbb{N}$ and let $C_{n, d}$ be defined as in \eqref{eq:cycle_graph}.
    Then, 
    \begin{equation*}
        a_1(C_{2n, d}) / a_1(C_{n, d}) \to  1/ 4 \quad \text{as} \quad n \to \infty.
    \end{equation*}
    \label{pro:bounds_ratio}
\end{pro}
\begin{pf}
    The claimed asymptotic behavior is derived using Lemma \ref{lem:cycle_spectrum} and the Taylor series expansion of the cosine function.
    As $n \to \infty$,
    \begin{equation*}
         \frac{a_1(C_{2n, d})}{a_1(C_{n, d})} = \frac{\sum_{k=1}^d 2\left(1 - \cos(\tfrac{k \pi}{n})\right)}{\sum_{k=1}^d 2\left(1 - \cos(\tfrac{2 k \pi}{n})\right)} \to \frac{\sum_{k=1}^d \left(\frac{k \pi}{2 n}\right)^2}{\sum_{k=1}^d \left(\frac{k \pi}{n}\right)^2} \to \frac{\frac{1}{4} \sum_{k=1}^d k^2}{\sum_{k=1}^d k^2} = \frac{1}{4}.
         \qed
    \end{equation*}
\end{pf}

To end this section, we provide Table \ref{tab:bounds_comparison}, that compares the bounds \eqref{eq:diameter_inequalities_2} and \eqref{eq:path_cycle} with existing ones.
This is done for $P_{n, d}$ with $d = 2$ (Table \ref{tab:bounds_comparison_a}) and $d = 3$ (Table \ref{tab:bounds_comparison_b}).
It can be noted that the proposed bounds outperform the existing ones for the graphs considered, with \eqref{eq:path_cycle} being the tightest.
Notably, \eqref{eq:diameter_inequalities_2}, although formulated for arbitrary graphs, provides a tighter bound for $P_{n, d}$ than the one given in \citet[Prop. 7.8]{Lew2025}.

\begin{table}[ht]

\centering
\begin{subtable}{\linewidth}
    \centering
    \begin{tabular}{c|ccccc}
    \toprule
    $n$ & $5$ & $15$ & $25$ & $35$ & $45$ \\
    \toprule
    $a_1(P_{n, 2})$  & 1.586 & 0.216 & 0.079 & 0.04 & 0.024 \\
    \midrule
    \citet[Th. 2.7]{Alon1985} & 43.131 & 9.968 & 4.792 & 2.913 & 1.994 \\
    \citet[Lem. 1.14]{Chung1996} & N/A & 0.667 & 0.445 & 0.354 & 0.304 \\
    New Bound \eqref{eq:diameter_inequalities_2} & \textbf{3.5} & \textbf{0.454} & \textbf{0.161} & \textbf{0.081} & \textbf{0.049}\\
    \midrule
    \citet[Prop. 7.8]{Lew2025} & 5.0 & 0.835 & 0.31 & 0.16 & 0.097 \\
    New Bound \eqref{eq:path_cycle} & \textbf{1.764} & \textbf{0.217} & \textbf{0.079} & \textbf{0.04} & \textbf{0.024} \\
    \bottomrule
    \end{tabular}
    \caption{$d = 2$}
    \label{tab:bounds_comparison_a}
\end{subtable} 
\begin{subtable}{\linewidth}
    \centering
    \begin{tabular}{c|ccccc}
    \toprule
    $n$ & $5$ & $15$ & $25$ & $35$ & $45$ \\
    \toprule
    $a_1(P_{n, 3})$  & 3.0 & 0.59 & 0.218 & 0.112 & 0.068 \\
    \midrule
    \citet[Th. 2.7]{Alon1985} & 43.131 & 29.306 & 16.174 & 8.77 & 6.434 \\
    \citet[Lem. 1.14]{Chung1996} & N/A & 0.953 & 0.691 & 0.546 & 0.487 \\
    New Bound \eqref{eq:diameter_inequalities_2} & \textbf{4.5} & \textbf{1.418} & \textbf{0.595} & \textbf{0.246} & \textbf{0.162}\\
    \midrule
    \citet[Prop. 7.8]{Lew2025} & 5.0 & 2.217 & 0.852 & 0.443 & 0.27 \\    
    New Bound \eqref{eq:path_cycle} & \textbf{4.382} & \textbf{0.599} & \textbf{0.219} & \textbf{0.112} & \textbf{0.068} \\
    \bottomrule
    \end{tabular}
    \caption{$d = 3$}
    \label{tab:bounds_comparison_b}
\end{subtable}
\caption{Comparison of upper bounds for the algebraic connectivity of $P_{n, d}$ for increasing $n$. Values are computed numerically and shown with $3$ decimal digits of precision.
The first three bounds of each (\subref{tab:bounds_comparison_a}) and (\subref{tab:bounds_comparison_b}) are applicable to arbitrary graphs, while the last two are specific for $P_{n, d}$.}
\label{tab:bounds_comparison}
\end{table}

\section{Concluding Remarks}
\label{sec:conclusion}
This paper presents new contributions to the generalized algebraic connectivity---a quantitative measure of graph rigidity in Euclidean spaces---and its interplay with fundamental graph invariants.
Recognizing the inherent connection between rigidity and connectivity, our investigation yields new metrics, bounds, and structural characterizations.
These results may be particularly useful in applications where both communication robustness (connectivity) and structural stability (rigidity) are important considerations.

The introduction of the $d$-rigidity ratio provides a novel perspective on quantifying a graph's rigidity.
By showing that $0 \leq \rho_d(G) \leq 1$, we offer a normalized measure that enables more direct comparisons of rigidity across different graphs.
We also examined how rigidity interacts with other structural properties of graphs, such as the diameter and the vertex connectivity, improving on previous characterizations and bounds. 
These results contribute to a better understanding of the spectral behavior of rigid networks.

Future research directions include studying the existence of constants $0 < \alpha_d, \beta_d < 1$ such that
\begin{equation*}
\alpha_d \leq \left\{ \rho_d(G) : G \text{ is rigid in } \mathbb{R}^d \right\} \leq \beta_d.    
\end{equation*}
Also, exploring the potential decreasing monotonicity of the $d$-rigidity ratio with respect to the dimension $d$, as conjectured by \citet[Conj. 8.3]{Lew2025}.
Finally, it would be of interest to investigate the computation and practical implications of the $d$-rigidity ratio in application domains such as network design, sensor placement, and multi-robot coordination.

\textcolor{black}{
\section*{Acknowledgments}
The authors thank the anonymous reviewer for their careful reading and insightful comments, which greatly contributed to improving the quality and clarity of this work.}

\bibliographystyle{elsarticle-num-names} 
\bibliography{cas-refs}






\end{document}